\documentclass{article}
\usepackage{amsrefs,amssymb,amsmath,amsthm,amsfonts,epsfig,graphicx}
\usepackage[utf8]{inputenc}
\usepackage[T1]{fontenc}
\usepackage[english]{babel}
\usepackage{hyperref}
\hypersetup{
    colorlinks,
    citecolor=black,
    filecolor=black,
    linkcolor=black,
    urlcolor=black
}

\theoremstyle{plain}
\newtheorem{teo}{Theorem}[section]

\newtheorem{lem}[teo]{Lemma}
\newtheorem{prop}[teo]{Proposition}
\newtheorem{rmk}[teo]{Remark}
\newtheorem{exa}[teo]{Example}

\theoremstyle{definition}
\newtheorem{df}[teo]{Definition}

\DeclareMathOperator{\per}{Per}

\DeclareMathOperator{\Emb}{Emb}

\DeclareMathOperator{\mesh}{mesh}

\DeclareMathOperator{\ord}{ord}

\DeclareMathOperator{\diam}{diam}
\DeclareMathOperator{\clos}{clos}

\DeclareMathOperator{\dist}{dist} 
\DeclareMathOperator{\distF}{d}

\newcommand{\aca}{$\bigstar$}
\newcommand{\U}{{\cal U}}
\newcommand{\V}{{\cal V}}

\newcommand{\R}    {\mathbb R}

\newcommand{\Z}  {\mathbb Z}
\newcommand{\N}  {\mathbb N}
\renewcommand{\epsilon}{\varepsilon}

\author{Alfonso Artigue}

\title{Partially Expansive Homeomorphisms}

\begin{document}
\date{\today}
\maketitle

\begin{abstract}
 In this article we consider several forms of expansivity. 
 We introduce two new definitions related with topological dimension. 
 We study the topology of local stable sets under cw-expansive surface homeomorphisms and 
 expansive homeomorphisms of three-manifolds. 
 The problem of $C^r$-stable expansivity is also considered.
\end{abstract}

%
%


\section{Introduction}

Given a homeomorphisms $f\colon X\to X$ of a compact metric space $(X,\dist)$ 
define the set $\Gamma_\delta(x)=\{y\in X:\dist(f^n(x),f^n(y))\leq\delta,\, \forall n\in\Z\}$, 
for $\delta>0$ and $x\in X$. 
This set is usually called as the \emph{dynamical ball}.
We say that $f$ is \emph{expansive} if for some $\delta>0$ it holds that 
$\Gamma_\delta(x)=\{x\}$ for all $x\in X$. 
Some variations of this definition consider special 
properties on the set $\Gamma_\delta(x)$. For example: 
\cite{Ka93} requires that $\Gamma_\delta(x)$ has vanishing topological dimension, 
\cite{Mo12} is related with the cardinality of $\Gamma_\delta(x)$, 
\cite{MoSi} considers that $\Gamma_\delta(x)$ has vanishing measure and 
\cite{Bo} asks this set to have vanishing topological entropy.


Let us explain another way of thinking about the definition of expansive homeomorphism 
that I learned from J. Lewowicz.
Consider that a $\delta$-observer is someone or something whose observations and measures 
have a precision of $\delta$. 
In this way, if two points have a distance smaller than $\delta$, such an observer will 
not be able to guess that these two points are different. 
Then, expansiveness, in terms of our $\delta$-observer, means that at some time he will \emph{see} 
two points, given that he knowns all of the iterates of this two different points. 
The idea of the $\delta$-observations was used in \cite{Ar14}, 
where it is asked to our $\delta$-observer to guess the cardinality of a finite set knowing its iterates. 
This led us in \cite{Ar14} to the definition of $(m,n)$-expansive homeomorphism 
that extends the concept of $N$-expansivity introduced by Morales \cite{Mo12}.

In this article we consider both ways of understanding the definition: 
in terms of the dynamical ball and with our $\delta$-observer. 
Both concepts are stated from the viewpoint of topological dimension.  
We consider that there is a natural number $D$ such that 
for all $x\in X$ we have that the topological dimension of $\Gamma_\delta(x)$ is less than $D$. 
This kind of homeomorphisms will be called partially expansive with central dimension $D$. 
We will also ask to our $\delta$-observer to guess the topological dimension of a set. 
The formalization of this idea will be called dimension-wise expansivity. 

This article also deals with standard versions of expansivity as 
the original one introduced by Utz and 
cw-expansivity defined by Kato, that in turn, are special instances of partial expansivity.
It is known that there are surface homeomorphisms that are cw-expansive but not expansive. 
We will show that this kind of homeomorphisms always presents what we call doubly asymptotic 
sectors. Expansive homeomorphisms are known to be conjugate with pseudo-Anosov diffeomorphisms. 
As we will see, the world of cw-expansive surface homeomorphisms with doubly asymptotic sectors 
has some anomalous examples. 
We will construct a cw-expansive surface homeomorphism with points whose local stable set is not a finite union of embedded arcs. 

For expansive homeomorphisms of three-manifolds, Vieitez \cite{Vi2002} obtained 
a conjugation with a linear Anosov diffeomorphism assuming some regularity and that there are no 
wandering points. 
In \cite{Vi2002} the following question is presented:
without assuming that $\Omega(f)=M$, can we have points with local stable and
unstable sets that are not manifolds for $f\colon M\to M$, an expansive diffeomorphism
defined on a three-dimensional manifold $M$?
We give a positive answer for homeomorphisms.

We also consider robust cw-expansiveness, i.e., diffeomorphisms in the $C^r$-interior of 
the set of $C^r$ cw-expansive diffeomorphisms of a compact surface. 
For $r=1$ Mañé \cite{Ma75} proved that a $C^1$ diffeomorphism of a compact manifold is 
robustly expansive if and only if it is quasi-Anosov (i.e., the tangent action is expansive). 
This result was extended by Sakai for cw-expansiveness \cite{Sa97} and measure expansiveness \cite{SaSuYa} (see also \cite{ArCa}). 
In this article we consider this problem in the $C^r$ topology with $r>1$. 
We introduce a generalization of quasi-Anosov diffeomorphisms allowing tangencies between 
stable and unstable arcs, but controlling the order of contact. 
In this way we find $C^r$ robustly cw-expansive diffeomorphisms not being quasi-anosov. 
These surface diffeomorphisms are in fact robustly $N$-expansive (in the sense of Morales).

Let us now describe the contents of this paper. 
In Section \ref{secTopDyn} we start recalling some basic properties of topological dimension. 
We introduce the definition of partially expansive homeomorphism with central dimension $D$. 
For $D=-1$ we prove the equivalence with expansivity (in the usual sense of Utz). 
For $D=0$ we show that partial expansiveness is equivalent with cw-expansiveness (in the sense of Kato). 
In the smooth category, it is known that the time-one diffeomorphism of a hyperbolic vector field is 
partially hyperbolic. In our topological setting we show that 
the time-one homeomorphisms of an expansive flow is partially expansive.
For a homeomorphism of a compact manifold of dimension $n$, we show that partial expansivity with central dimension $D=n-1$ 
is equivalent with what is called as sensitivity to initial conditions. 
Next we consider another definition that we call dw-expansivity. 
We show that this definition is stronger than partial expansivity. 
We show that the converse does not hold via our previous result related with expansive flows. 
In the case of vanishing central dimension we prove that both definitions are equivalent with cw-expansivity. 
We also consider positive dw-expansive endomorphisms. We show that every positive expansive endomorphism 
is positive dw-expansive. Via inverse limits we obtain examples of positive dw-expansive homeomorphisms. 

In Section \ref{secCwexp} we consider cw-expansive surface homeomorphisms. The main result of this section states that 
every cw-expansive surface homeomorphism without doubly asymptotic sectors is in fact expansive. 
This result is based on the work of Lewowicz \cite{L}. 
In this work he proved that expansive surface homeomorphisms are conjugate with pseudo-Anosov diffeomorphisms. 
In the first part of his article, it is constructed a local product structure around every point excluding a finite set. 
We noticed that the hypothesis of expansivity is used in two ways (in \cite{L}). 
First it is used to prove that there are no Lyapunov stable points and then to construct continua of uniform size 
contained in the local stable set of every point. 
After the results of Kato on cw-expansive homeomorphisms we know that this can be done only assuming cw-expansivity. 
Then, in \cite{L}, expansivity is used to avoid that a local stable continuum meets twice a local unstable continuum. 
This is what we call a doubly asymptotic sector. In light of this remarks, Theorem \ref{TeoEquivExp} is natural. 
The techniques of this section are essentially contained in \cite{L}, but a reorganization of the results allows us 
to easily exclude the existence of spines.

In Section \ref{secAno} we present examples that limit the extension of the results of the previous section. 
An important step in the construction of a local product structure, as before, is to prove that 
local stable continua are locally connected. 
We give an example of a cw-expansive surface homeomorphism with a point 
whose local stable set is connected but not locally connected. 
Of course, it has doubly asymptotic sectors. 
We also construct an expansive homeomorphism on a three-manifold with similar properties. 
In both cases there are wandering points and the ideas are based on the quasi-Anosov diffeomorphism 
that is not Anosov defined in \cite{FR}.

In Section \ref{secSta} we start recalling known facts related with Axiom A diffeomorphisms 
and with what we call $\Omega$-expansivity. 
We introduce the definition of $Q^r$-Anosov $C^r$-diffeomorphisms extending the concept of quasi-Anosov 
diffeomorphism. We prove that $Q^r$-diffeomorphisms are $C^r$-robustly $r$-expansive ($N$-expansive in the sense of Morales with $N=r$). 
In particular they are cw-expansive.
We give an example showing that such diffeomorphisms do exist.

\section{Topological dynamics}
\label{secTopDyn}
In this section we will introduce two concepts for the dynamics of homeomorphisms on compact metric spaces: partial expansiveness and dimension expansiveness. 
The former satisfies that the time-one homeomorphism of an expansive flow is a partially expansive 
homeomorphism. The later is related with what can be called the $\delta$-\emph{observable dimension}.
Both definitions are stated in terms of topological dimension, therefore we will start recalling 
the definitions and known results. The interested reader should consult \cite{HW} for more on this subject.

\subsection{Topological dimension}

Let $(X,\dist)$ be a compact metric space. 
Denote by $\tau$ the topology of $X$, i.e. the set of open subsets of $X$.  
For $\U\subset \tau$ we will define the \emph{order} of $\U$.
We say that
$\ord \U\leq n$ if for all 
${\cal V}\subset\U$ with $|{\cal V}|>n+1$ it holds that $\cap {\cal V}=\emptyset$.
If $\ord\U\leq n$ and it is not true that $\ord\U\leq n-1$ we say that $\ord\U=n$. 
Given $\U\subset\tau$ we define 
$$
  \mesh (\U)=\sup_{U\in\U} \diam(U).
$$
For $\epsilon>0$ and a compact set $Y\subset X$ define the $\epsilon$-\emph{dimension}
$$
\dim_\epsilon Y=\inf\{\ord\U:\U\subset \tau, Y\subset\cup\U,\mesh(\U)<\epsilon\}.
$$

\begin{rmk}
The $\epsilon$-dimension of $Y\subset X$  can also be defined considering the relative topology 
of $Y$ as a topological subspace of $X$. 
Let us prove that these two numbers are equal. 
Given a compact set $Y$ and two open sets $U,V$ such that $U\cap Y\cap V=\emptyset$ we have to show that 
there are two open sets $U', V'$ such that $U\cap Y=U'\cap Y$, $V\cap Y=V'\cap Y$ and 
$U'\cap V'=\emptyset$. Define 
$$U'=\{x\in U:\dist(x,Y\cap U)<\dist(x,Y\cap V)\}$$
and
$$V'=\{x\in V:\dist(x,Y\cap V)<\dist(x,Y\cap U)\}.$$
\end{rmk}

Note that if $\epsilon_1<\epsilon_2$ then $\dim_{\epsilon_1} Y\geq\dim_{\epsilon_2} Y$, 
therefore $\lim_{\epsilon\to 0}\dim_\epsilon(Y)$ is an integer number or $\infty$. 
Define the \emph{(topological) dimension} of $Y$ as 
$$\dim Y=\lim_{\epsilon\to 0}\dim_\epsilon(Y).$$
By the \emph{Covering Theorem for Compact Spaces} in Chapter V, Section 8 of \cite{HW}, it is the 
usual definition of topological dimension. 

\begin{prop}
\label{dimeps}
 For every compact set $Y\subset X$ with $\dim(Y)<\infty$ there is $\epsilon>0$ such that 
 $\dim_\epsilon(Y)=\dim(Y)$.
\end{prop}

\begin{proof}
 It follows by definition because $\dim_\epsilon(Y)$ and $\dim(Y)$ are integer numbers.
\end{proof}

In the plane we can consider a sequence of arcs converging, in the Hausdorff metric, to a point. 
Also there are sequences of arcs converging to a disk. Therefore the operator $\dim$ is not continuous. 
For $\dim_\epsilon$ we have the following result.

\begin{prop}
\label{dimlim}
 If $C_n$ is a sequence of compact subsets of $X$ and $C_n\to C$ in the Hausdorff metric then for every fixed $\epsilon>0$ we have that
 $$\limsup_{n\to\infty} \dim_\epsilon(C_n)\leq\dim_\epsilon(C).$$
\end{prop}

\begin{proof}
 Let $d=\dim_\epsilon(C)$. By definition there is $\U\subset \tau$ such that $C\subset\cup\U$, $\mesh(\U)<\epsilon$ and 
 if $\V\subset \U$ and $|\V|>d+1$ then $\cap \V=\emptyset$. 
 Take $n_0$ such that $C_n\subset \cup \V$ for all $n\geq n_0$. 
 Then $\dim_\epsilon(C_n)\leq d$ for all $n\geq n_0$ and the proof ends.
\end{proof}

Notice that every compact set $C$ can be approximated by a sequence $C_n$ of finite subsets of $C$. So, it could be the case that
$\limsup_{n\to\infty} \dim_\epsilon(C_n)<\dim_\epsilon(C).$

\begin{prop}
\label{dimcc}
 For every compact metric space $Y$ it holds that 
 \[
  \dim(Y)=\sup\{\dim(C):C\subset Y, C\hbox{ is a continuum}\}.
 \]
\end{prop}

\begin{proof}
 Let $L=\sup\{\dim(C):C\subset Y, C\hbox{ is a continuum}\}$.
 In general, if $A\subset B$ it holds that $\dim(A)\leq\dim(B)$. 
 If $L=\infty$ the result is trivial. Then we will assume that $L<\infty$.
 It is sufficient to prove that there is a component $C$ of $Y$ with $\dim(C)=\dim(Y)$. 

 By contradiction assume that every component of $Y$ has dimension smaller than the dimension of $Y$. 
 Take $\epsilon >0$. 
 Given $x\in Y$ consider the component $C_x$ containing $x$. 
 Consider $\U_x\subset \tau$ such that 
 $\mesh(\U_x)<\epsilon$, $\ord (\U_x)<\dim (Y)$ and $C_x\subset\cup \U_x$. 
 We can also assume that $K_x=\cup\U_x$ is compact (see Theorem 2-15 in \cite{HY}). 
 Take $x_1,\dots,x_n$ such that $Y=\bigcup_{i=1}^n \cup \U_{x_i}$. 
Define $\V_1=\U_{x_1}$, $\V_2=\{U\setminus K_1:U\in \U_{x_2}\}$ and in general 
\[
 \V_i=\{U\setminus(K_1\cup\dots\cup K_{i-1}):U\in\U_{x_i}\}
\]
for $i=2,\dots,n$. 
Define $\V=\cup_{i=1}^n\V_i$. 
Now it is easy to see that $\ord(\V)\leq L<\dim(Y)$, $\cup \V=Y$ and $\mesh(\V)<\epsilon$. 
This is a contradiction because $\epsilon$ is arbitrary. 
\end{proof}

\subsection{Partial expansiveness}

Let $(X,\dist)$ be a compact metric space. 
We say that $C\subset X$ is \emph{non-trivial} if it has 
more than one point.
If $C$ has just one point we say that $C$ is a \emph{singleton}.

\begin{df}
 Given an integer $D\geq -1$, we say that a homeomorphism $f\colon X\to X$ is \emph{partially expansive} 
 with \emph{central dimension} $D$ and \emph{expansive constant}
 $\epsilon>0$ 
 if for every non-trivial compact set $C\subset X$ with $\dim(C)>D$ 
 there is $k\in\Z$ such that 
  $\diam(f^k(C))\geq\epsilon$.
\end{df}

Let us give some examples. 
Recall that a homeomorphism $f\colon X\to X$ is \emph{expansive} 
with \emph{expansive constant} $\epsilon>0$ if 
$\dist(f^n(x),f^n(y))<\epsilon$ for all $n\in\Z$ implies $x=y$. 

\begin{prop}
A homeomorphism $f\colon X\to X$ is expansive if and only if it is 
partially expansive with central dimension $D=-1$.
\end{prop}

\begin{proof}
Let us start proving the direct part. 
Assume that $f$ is expansive with expansive constant $\epsilon$. 
Let $C\subset X$ be a non-trivial compact set. 
Since $C$ is non-trivial we have that $\dim(C)\geq 0$ and there are 
$x,y\in C$ with $x\neq y$. 
Since $f$ is expansive there is $k\in\Z$ such that $\dist(f^k(x),f^k(y))\geq\epsilon$. 
Therefore $\diam(f^k(C))\geq\epsilon$. So, $f$ is partially expansive with central dimension $D=-1$.

 In order to prove the converse assume that 
 $f$ is partially expansive with central dimension $D=-1$. Take $\epsilon>0$ an expansive constant. 
 Given $x\neq y$ consider the compact set $C=\{x,y\}$.
 Since $\dim(C)=0>-1$ 
 there is $k\in\Z$ so that $\diam(f^k(C))\geq\epsilon$. 
 Therefore $\dist(f^k(x),f^k(y))\geq\epsilon$ and $f$ is expansive. 
\end{proof}

A homeomorphism $f\colon X\to X$ is \emph{cw-expansive} if there is $\delta>0$ such that 
if $C\subset X$ is a continuum and
 $\diam(f^j(C))<\delta$ for all $j\in\Z$ then $C$ is a singleton.

\begin{prop}
\label{cwPexpD0}
A homeomorphism $f\colon X\to X$ is cw-expansive if and only if it is partially expansive with central dimension $D=0$.
\end{prop}

\begin{proof}
 It follows by Proposition \ref{dimcc} above and Proposition 2.3 in \cite{Ka93}.
\end{proof}

Let $\phi\colon\R\times X\to X$ be a continuous flow. 

\begin{df}
We say that a flow $\phi$ is \emph{kinematic expansive} if for all $\epsilon>0$ there is $\delta>0$ such that if 
$\dist(\phi_t(x),\phi_t(y))<\delta$ for all $t\in\R$ then $y=\phi_s(x)$ for some $s\in(-\epsilon,\epsilon)$. 
\end{df}

Examples of kinematic expansive flows are expansive flows in the sense of Bowen-Walters \cite{BW} and $k^*$-expansive flows as defined by Komuro \cite{Ko84}. 
See \cite{ArKinExp} for more on kinematic expansive flows.

\begin{prop}
\label{propKexp}
If $\phi$ is a kinematic expansive flow then
for all $T\in\R$, $T\neq 0$, the homeomorphism 
$\phi_T\colon X\to X$ is partially expansive with central dimension $D=1$.
\end{prop}

\begin{proof}
Fix $T\neq 0$. 
Without loss of generality assume that $T>0$.
Let $\sigma>0$ be such that if $\dist(\phi_s(x),\phi_s(y))<\sigma$ for all $s\in\R$ 
then $y\in \phi_\R(x)$. 
Let $\delta>0$ be such that 
if $C\subset X$ is a continuum and $\diam(C)<\delta$ then $\diam(\phi_s(C))<\sigma$ 
for all $s\in[-T,T]$. 
Then, if $\diam(\phi_{nT}(C))<\delta$ for all $n\in\Z$ we have that 
$\diam(\phi_s(C))<\sigma$ for all $s\in\R$. 
Therefore, $C$ is an orbit segment and consequently $\dim(C)\leq 1$. 
We have proved that $\delta$ is an expansive constant for $f=\phi_T$ with central dimension 
$D=1$.
\end{proof}

Recall that $\phi$ is an \emph{expansive flow} (in the sense of Bowen and Walters) if for all 
$\epsilon>0$ there is $\delta>0$ such that if $\dist(\phi_t(x),\phi_{h(t)}(y))<\delta$ 
for all $t\in\R$, with $h\colon\R\to\R$ an increasing homeomorphism satisfying $h(0)=0$, 
then there is $s\in\R$ such that $y=\phi_s(x)$ and $|s|<\epsilon$. 
Also recall that $\phi$ is $k^*$-\emph{expansive} if for all $\epsilon>0$ there is 
$\delta>0$ such that if $\dist(\phi_t(x),\phi_{h(t)}(y))<\delta$ 
for all $t\in\R$, with $h\colon\R\to\R$ an increasing homeomorphism satisfying $h(0)=0$, 
then there $x$ and $y$ are contained in an orbit segment of diameter smaller than $\epsilon$.
Expansive flows and $k^*$-expansive flows 
are kinematic expansive. 
Therefore $\phi_T$ is a partially expansive homeomorphism 
with $D=1$ if $\phi$ is expansive or $k^*$-expansive.

\begin{df}
 We say that a homeomorphism $f\colon X\to X$ 
 is \emph{sensitive to initial conditions}
 if there is $\rho>0$ such that for all $x\in M$ and for all $r>0$ there are $y\in B_r(x)$ and 
  $n\in\Z$ such that $\dist(f^n(y),f^n(x))>\rho$.
\end{df}

\begin{prop}
\label{sensitive}
Let $f\colon M\to M$ be a homeomorphism of a compact manifold $M$ of dimension $n$. 
 The following statements are equivalent:
 \begin{enumerate}
  \item $f$ is partially expansive with central dimension $D=n-1$,
  \item $f$ is sensitive to initial conditions.
 \end{enumerate}
\end{prop}

\begin{proof}
$(1\Rightarrow 2)$. Let $\delta$ be an expansive constant and define $\rho=\delta/2$. 
Given $x\in M$ and $r>0$ define $C=\clos(B_r(x))$. 
Since $M$ has dimension $n$ we have that $\dim(C)=n$ (see Corollary 1, page 46 in \cite{HW}). 
Then there is $n\in\Z$ such that $\diam(f^n(C))>\delta$. 
Therefore, there is $y\in B_r(x)$ such that $\dist(f^n(x),f^n(y))>\rho$.

$(2\Rightarrow 1)$. 
Suppose now that $f$ is sensitive to initial condition and take $\rho$ from the definition. 
Let $C\subset X$ be a compact subset such that $\dim(C)=n$. 
By the cited result from \cite{HW} we have that $C$ has non-empty interior. 
Take $x$ in the interior of $C$ and $r>0$ such that 
$B_r(x)\subset C$. 
We know that there are $y\in B_r(x)$ and $n\in\Z$ such that 
$\dist(f^n(y),f^n(x))>\rho$. 
Therefore, $\diam(f^n(C))>\rho$, which proves that $\rho$ is an expansive constant.
\end{proof}

\subsection{Dimension-wise expansive homeomorphisms}

In this section we introduce \emph{dw-expansiveness} as a stronger form of partial expansiveness. 

\begin{df}
 We say that $f\colon X\to X$ is \emph{dw-expansive} 
 if there are $D\geq 0$ and $\sigma>0$ such that if 
 $\dim(C)>D$ for a compact set $C\subset X$ then there is 
 $n\in\Z$ such that $\dim_\sigma(f^n(C))>D$.
 In this case we say that $D$ is the \emph{central dimension} and $\sigma$ is a \emph{dw-expansive constant}.
\end{df}

Dw-expansiveness means that given a set with dimension greater than $D$ the $\delta$-observer will 
see, at some iterate of the set, that it has dimension $D$, i.e. $\dim_\delta(f^n(D))>D$.
Notice that he will say that the 
set $[0,\epsilon]\times [0,7]\subset \R^2$ is one-dimensional if $\epsilon\in(0,\delta)$. 
\begin{prop}
\label{propDwPexp}
Every dw-expansive homeomorphism is partially expansive with the same central dimension.
\end{prop}

\begin{proof}
 Let $f$ be a dw-expansive homeomorphism and let $D,\sigma$ be as in the definition. 
 Notice that if $\dim_\sigma(C)>D$ then $\diam(C)>\sigma$. 
 Therefore, $\sigma$ is an expansive constant for $f$ with central dimension $D$.
\end{proof}

\begin{exa}
Let us prove that the converse of Proposition \ref{propDwPexp} does not hold. 
Consider the annulus 
\[
 A=\{(x,y)\in\R^2: 1\leq \sqrt{x^2+y^2}\leq 2
\]
and the vector field $X\colon\R^2\to\R^2$ given by 
\[
 X(x,y)=\frac{(-y,x)}{\sqrt{x^2+y^2}}.
\]
Let $\phi\colon\R\times A\to A$ be the flow of $X$.  
Every orbit is periodic and the period increases with the distance to the origin. 
Therefore, $\phi$ is kinematic expansive.
Define $f=\phi_1$ the time-one diffeomorphism of 
the annulus. 
By Proposition \ref{propKexp} we have that $f$ is partially expansive with central dimension $D=1$. 
Let us show that $f$ is not dw-expansive with central dimension $D=1$. 
Given $\sigma>0$ consider 
\[
 C=\{(x,y)\in\R^2: 1\leq \sqrt{x^2+y^2}\leq 1+\epsilon\}
\]
for some $\epsilon>0$. 
We have that $\dim(C)=2$, $f^n(C)=C$ for all $n\in\Z$ and 
if $\epsilon $ is small it is easy to see that $\dim_\sigma(C)=1$. 
Therefore, $f$ is not dw-expansive with central dimension $D=1$.
\end{exa}

\begin{prop}
 The following statements are equivalent:
 \begin{enumerate}
  \item $f$ is dw-expansive with central dimension $D=0$
  \item $f$ is partially expansive with central dimension $D=0$
  \item $f$ is cw-expansive
 \end{enumerate}
\end{prop}

\begin{proof}
 ($1\Rightarrow 2$) follows by Proposition \ref{propDwPexp} and 
 ($2\Leftrightarrow 3$) is proved in Proposition \ref{cwPexpD0}. 
 Now we will prove that $(3\Rightarrow 1)$.
 Let $f\colon X\to X$ be a cw-expansive homeomorphism with cw-expansive constant $\delta$. 
 Take a compact set $C\subset X$ such that $\dim(C)>0$. 
 Let $A\subset C$ be a non-trivial connected component of $C$. 
 Since $f$ is cw-expansive there is $n\in\Z$ such that 
 $\diam(f^n(A))>\delta$. Then $\dim_\delta(A)\geq\dim_\delta(C)\geq 1$. 
 Therefore, $\delta$ is a dw-expansive constant.
\end{proof}

\subsubsection{Positive dw-expansivity}

Let $(X,\dist)$ be a compact metric space. 
In this section we will consider a continuous map $f\colon X\to X$ that is an open map and $f(X)=X$. 

\begin{df}
 We say that a map $f\colon X\to X$ is \emph{positive dw-expansive} 
 with \emph{central dimension} $D$ if there is $\sigma>0$ 
 such that if $\dim(C)>D$ then there is $n\geq 0$ such that
 $\dim_\sigma(C)>D$.
\end{df}

Recall that $f\colon X\to X$ is \emph{positive expansive} 
if there is $c>0$ such that if $x\neq y$ then there is $n\geq 0$ such that $\dist(f^n(x),f^n(y))>c$. 

\begin{rmk}[Adapted metric] 
\label{RmkAdapted}
Since we are assuming that $f$ is an open map, we have by \cite{Reddy} 
that every positive expansive map $f$ is Ruelle-expanding, that is:  
there are a metric $\distF$ defining the topology of $X$ and
constants $r>0$ and
$\lambda>1$ such that for every $x\in X$ and $a\in f^{-1}(x)$ there  
exists $\varphi\colon B_r(x)\to X$ such that $\varphi (x)=a$, $f\circ \varphi(y)=y$ 
for every $y\in B_r(x)$ and $\lambda\distF (\varphi (z),\varphi(w))\leq \distF (z,w)$, 
for every $z, w\in B_r(x)$.
We say that $\distF$ is an \emph{adapted metric}.
\end{rmk}

\begin{prop}
\label{posdw}
 Every positive expansive map is positive 
 dw-expansive with any central dimension 
 $D\geq 0$.
\end{prop}

\begin{proof}
As in Remark \ref{RmkAdapted} consider an adapted metric $\distF$ and the constants $c>0$ and $\lambda>1$. 
Let $\U$ be a finite covering of $X$ with $\mesh(\U)<r/2$.
Let $\sigma>0$ be a Lebesgue number of $\U$, i.e., if $A\subset X$ and $\diam(A)$ then there is 
$U\in\U$ such that $A\subset U$.
Fix $D\geq 0$ and a compact set $C\subset X$ such that 
$\dim_\sigma(f^n(C))\leq D$ for all $n\geq 0$. 
We will show that $\dim(C)\leq D$. 
For this, consider $\epsilon>0$ and take $m\geq 0$ such that 
$r/\lambda^m<\epsilon$. 
Since $\dim_\sigma(f^m(C))\leq D$ there is a
covering $\U_1$ of $f^m(C)$ with $\mesh(\U_1)<\sigma$ and $\ord(\U_1)\leq D$. 
Suppose that $\U_1=\{V_1,\dots,V_k\}$.
By Remark \ref{RmkAdapted} each $V_i\in\U_1$ we have that 
$f^{-m}(V_i)=V_i^1\cup\dots\cup V^{k_i}_i$ a disjoint union of open sets such that 
$\diam(V_i^j)<r/\lambda^m$ for all $j=1,\dots,k_i$. 
Define a covering of $C$ by 
\[
 \U_C=\{V_i^j:i=1,\dots,k;j=1,\dots,k_i\}.
\]
By construction we have that $\mesh(\U_C)<\epsilon$ and $\ord(\U_C)\leq D$. 
Since $\epsilon$ is arbitrary we conclude that $\dim(C)\leq D$. 
Therefore $f$ is positive dw-expansive with dw-expansive constant $\sigma$ and arbitrary central dimension $D\geq 0$.
\end{proof}

In light of the previous result it is natural to ask if every expansive homeomorphism is dw-expansive with arbitrary central dimension. 
Let us explain that this is not the case.

\begin{rmk}
 Let us give an example of an expansive homeomorphism that is not dw-expansive with central dimension $D=1$. 
 Consider the Anosov diffeomorphism $f$ of the two-dimensional torus given by $f(x,y)=(2x+y,x+y)$. 
 Given $\sigma>0$ consider a rectangle $C$ bounded by two small stable arcs and two small unstable arcs such that 
 $\diam(C)<\sigma$. It is easy to see that $\dim_\sigma(f^n(C))\leq 1$. Since $\dim(C)=2$ and $\sigma$ is arbitrary 
 we have that $f$ is not dw-expansive with central dimension $D=1$. 
 Notice that this argument holds for every pseudo-Anosov surface diffeomorphism.
\end{rmk}

\begin{rmk}[Inverse limits]
Given a continuous map $f\colon X\to X$ define $Y=\{a\in M^{\Z}: g(a_n)=a_{n+1}\}$ and $g\colon Y\to Y$ 
as $g(a)_n=a_{n+1}$ (the shift map). 
The distance in $Y$ is 
\[
\dist(a,b)=\sum_{n\in\Z}\frac{\dist(a_n,b_n)}{2^{|n|}}. 
\]
It is well known that $g$ is a homeomorphism, and it is called the \emph{inverse limit} of $f$.
With Proposition \ref{posdw} it is easy to see that the inverse limit of a positive expansive map is a positive 
dw-expansive homeomorphism with any central dimension 
$D\geq 0$.
\end{rmk}

Let us make the following question.
If a compact metric space $X$ admits a dw-expansive homeomorphism, is it true that 
$\dim(X)<\infty$?

\subsection{Plaque expansivity}

Let $f\colon X\to X$ be a homeomorphism on a compact metric space $(X,\dist)$. 
Given an equivalence relation $R$ on $X$ we denote by $R(x)$ the equivalence class of $x\in X$. 
Assume that $R(f(x))=f(R(x))$ for all $x\in X$.
Denote by $cc_x(A)$ the connected component of a set $A\subset X$ that contains $x$. 
For $\epsilon>0$ and $x\in X$ define $R_\epsilon(x)=cc_x(B_\epsilon(x)\cap R(x))$. 
An $(\epsilon,R)$-\emph{orbit} is a sequence $x_n\in X$, $n\in \Z$, such that 
$f(x_n)\in R_\epsilon(x_{n+1})$ for all $n\in\Z$.

\begin{df}
 Given a homeomorphism $f\colon X\to X$ of a compact metric space $(X,\dist)$ and an equivalence relation $R$ on $X$ 
 we say that $f$ is \emph{expansive mod} $R$ if for all $\epsilon>0$ there is $\delta>0$ such that 
 if $x_n, y_n$ are $(\epsilon,R)$-orbits and $\dist(x_n,y_n)<\delta$ for all $n\in \Z$ then 
 $x_n\in R_\epsilon(y_n)$ for all $n\in \Z$.
\end{df}

\begin{prop}
\label{propPexp}
If $f$ is expansive mod $R$ and for some $\epsilon>0$ it holds that 
$\dim(R_\epsilon(x))\leq D$ for all $x\in X$
 then $f$ is partially expansive with central dimension $D$.
\end{prop}

\begin{proof}
By definition there is $\delta$ corresponding to $\epsilon$. 
Suppose that $\diam(f^n(C))<\delta$ for all $n\in\Z$. 
Then, $C\subset R_\epsilon(x)$ for all $x\in C$. 
Consequently, $\dim(C)\leq D$. 
This proves that $f$ is partially expansive with central dimension $D$.
\end{proof}

\begin{rmk}
\label{obsPlExp}
Given a foliation $F$ of a compact manifold $M$ consider the equivalence relation of being in a common leaf. 
Define $F(x)$ as the leaf of $x$ (the equivalence class of $x$) and assume that $F$ is invariant by a homeomorphism $f\colon M\to M$.
Then $f$ is expansive mod $F$ if and only if $f$ is plaque expansive in the sense of \cite{HPS}.
In our topological setting $R_\epsilon(x)$ are the \emph{plaques}.
\end{rmk}

A $C^1$ diffeomorphism $f\colon M\to M$ is \emph{normally hyperbolic} 
to a $C^1$ foliation $F$ of $M$ if $f$ preserves $F$ and $TM=E^s\oplus TF\oplus E^u$, 
$df$ preserves the sub-bundles $E^s$ and $E^u$, 
and for a suitable Finsler on $TM$ it holds that
$\|df^{-1}_x|E^u\|<1$ and $\|df_x|E^s\|<1$ for all $x\in M$.

\begin{prop}
 If $f\colon M\to M$ is normally hyperbolic to a foliation $F$ and 
 $D$ is the dimension of the leaves of $F$ then $f$ is partially expansive with central dimension $D$.
\end{prop}

\begin{proof}
 By \cite[Theorem 7.2]{HPS} we know that $f$ is plaque expansive. 
 By Remark \ref{obsPlExp} we have that $f$ is expansive mod $F$. 
 Since $\dim(F_\epsilon(x))=D$ for all $\epsilon>0$ and for all $x\in M$ we have by Proposition \ref{propPexp} 
 that $f$ is partially expansive with central dimension $D$.
\end{proof}

A diffeomorphism $f\colon M\to M$ is \emph{partially hyperbolic} 
if it admits a $df$-invariant splitting of $TM=E^s\oplus E^c\oplus E^u$, 
for a suitable Finsler on $TM$ it holds that
$\|df^{-1}_x|E^u\|<1$ and $\|df_x|E^s\|<1$ for all $x\in M$ and 
$\|df_x v^s\|<\|df_x v^c\|<\|df_x v^u\|$ 
for all unit vectors $v^i\in E^i$, with $i\in \{s,c,u\}$.

\begin{rmk} 
If $f$ is a partially hyperbolic diffeomorphism 
with a $C^1$ foliation $F$ tangent to the central sub-bundle $E^c$
then $f$ is normally hyperbolic to $F$ and consequently, $f$ is partially expansive. 
\end{rmk}

\begin{rmk}
Given a continuous flow $\phi\colon\R\times X\to X$ we define the equivalence relation $O$ by $y\in O(x)$ 
if and only if there is $t\in\R$ so that $y=\phi_t(x)$, that is, $O(x)$ is the orbit of $x$. 
The following statements are equivalent: 
\begin{enumerate}
 \item $\phi_t$ is expansive mod $O$ for all $t\neq 0$,
 \item $\phi$ is expansive (in the sense of Bowen and Walters \cite{BW}),
\end{enumerate}
A proof can follow the ideas in \cite[Theorem 3]{BW}. Therefore, Bowen-Walters expansivity is 
a form of plaque expansivity, as noted in \cite[page 184]{BW}.
\end{rmk}

\section{Continuum-expansiveness on surfaces}
\label{secCwexp}
In this section we will show that for surface homeomorphisms expansiveness is equivalent with 
cw-expansiveness plus the non-existence of doubly asymptotic sectors. 

\subsection{Doubly asymptotic sectors and regularity}

Let $S$ be a compact surface without boundary and 
consider a cw-expansive homeomorphism $f\colon S\to S$ with expansive constant $\delta$. 
A continuum $C\subset S$ is \emph{stable} if $\diam(f^n(C))\leq\delta$ for all $n\geq 0$. 
In this case it holds that $\diam(f^n(C))\to 0$ as $n\to+\infty$. 
We say that $C$ is \emph{unstable} if it is stable for $f^{-1}$. 

\begin{df}
 We say that $f$ is a \emph{regular cw-expansive} homeomorphism 
 if every stable and every unstable continuum is arc-connected.
\end{df}

It is known that expansive surface homeomorphisms are regular (see \cites{Hi,L}). 
In Section \ref{AnomCwexp} we will show that there are non-regular 
cw-expansive surface homeomorphisms. 
If $f$ is expansive with expansive constant $\delta$, $C$ is stable and $D$ is unstable then $C\cap D$ has at most one point. 
This fact induces the following definition. 
We assume that $\delta$ is so small that every set of diameter 
smaller than $2\delta$ is contained in an open set $U\subset S$ with 
$U$ homeomorphic to $\R^2$.

\begin{df}
\label{defDas}
Let $C$ be a stable continuum and let $D$ be an unstable continuum. 
We say that $(C,D)$ form a \emph{(topological)  doubly asymptotic sector} 
if $C\cap D$ is not a singleton.
\end{df}

\begin{prop}
\label{cwNotbimpreg}
If $f\colon S\to S$ is cw-expansive and has no doubly asymptotic sectors then $f$ is regular.
\end{prop}

For the proof we follow Lemma 2.3 in \cite{L}, more details can be found there, we only give a sketch. 
See also Proposition 3.1 in \cite{Hi}.

\begin{proof}
 By contradiction assume that $A$ is a stable continuum that is not locally connected. 
 Let $A_n$, $n\geq 1$, be a sequence of pairwise disjoint stable sub-continua of $A$. 
 Assume that $A_n$ converges to a sub-continuum $A_*\subset A$ in the Hausdorff metric and $A_n$ is disjoint from $A_*$ for all 
 $n\geq 1$. 
 Take $p\in A_*$ and $\delta>0$ such that 
 $A_*$ and $A_n$ separate $B_\delta(p)$ for all $n\geq n_0$ (for some $n_0$ large). 
 See Figure \ref{figConLoc}.
 \begin{figure}[ht]
 \begin{center}
  \includegraphics{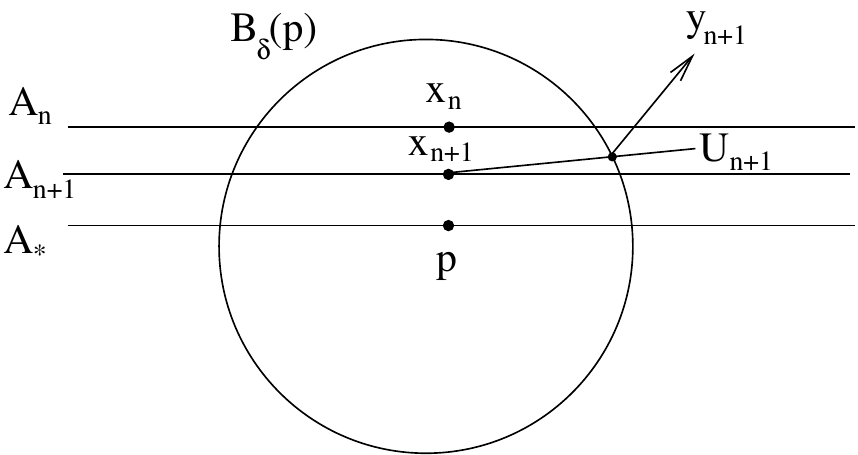}
  \caption{Local connection of stable sets.}
  \label{figConLoc}
 \end{center}
 \end{figure}
 Take $x_n\in A_n$ such that $x_n\to p$. 
 For each $n$ consider an unstable continuum $U_n$ not contained in $B_\delta(p)$ such that 
 $x_n\in U_n$. Since there are no doubly asymptotic sectors, each $U_n$ can 
 meet just one of the stable continua of the sequence $A_n$. 
 Take $y_n\in A_n\cap \partial B_\delta(p)$. 
 Let $U$ be limit continuum of $U_n$. 
 If $q$ is a limit point of $y_n$ we have that $\{p,q\}\subset U\cap A_*$. 
 Then, $U$ and $A_*$ make a doubly asymptotic sector. 
 This contradiction proves the result.
\end{proof}

\begin{rmk}
 As in \cites{Hi,L} a detailed proof of the previous result should use Janiszewski's Theorem on continua separating the plane.
\end{rmk}

\begin{prop}
\label{inestsepara}
Let $f\colon S\to S$ be a regular cw-expansive homeomorphism.
Suppose that there is a small disc $U\subset S$ without doubly asymptotic sectors 
and containing 
an unstable arc $\alpha$
in its boundary.
Then for all $x\in \alpha$ there is a stable arc from $x$ to $\partial U$ contained in $U$. 
\end{prop}

This proof is based in the proofs of Lemmas 3.1 and 3.2 in \cite{L}. Again, we only give the main ideas, the details can be found in \cite{L}.

\begin{proof} 
Given $x\in \alpha$ denote by $C$ the 
union of all stable continua contained in $\clos(U)$ and containing $x$.
Let $\beta$ be the arc $\partial U\setminus \alpha$.
By contradiction assume that $C\cap \beta=\emptyset$. 
The point $x$ separates $\alpha$ in two curves $\alpha_1$ and $\alpha_2$.
Given $\epsilon>0$ consider a curve $J$ starting at $\alpha_1$ and ending in $\alpha_2$ such that 
the interior of $J$ is contained in $(U\setminus C)\cap B_\epsilon(C)$. 
See Figure \ref{figSep}.
 \begin{figure}[ht]
 \begin{center}
  \includegraphics{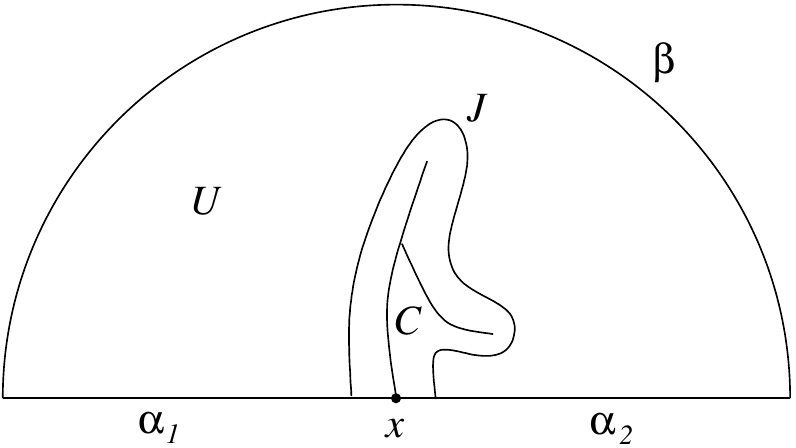}
  \caption{Illustration for the proof of Proposition \ref{inestsepara}.}
  \label{figSep}
 \end{center}
 \end{figure}
If $\epsilon$ is small we known that there is no stable arc from a point of $J$ to $\beta$. 
Therefore, for every $z\in J$ a stable continuum through $z$ meets the boundary of $U$ at $\alpha$. 
Define $J_i$, for $i=1,2$, as the set of points $y\in J$ for which there is a stable continuum through $y$ that cuts $\alpha_i$. 
Both sets, $J_1$ and $J_2$, are closed and non-empty. 
Also $J=J_1\cup J_2$. 
Since $J$ is connected, there is $p\in J_1\cap J_2$. 
This contradicts that there are no doubly asymptotic sectors in $U$.
\end{proof}

\subsection{Product boxes}

\begin{df}
A \emph{product box} is a homeomorphism $\phi\colon [0,1]\times[0,1]\to K\subset S$ 
taking horizontal lines onto stable arcs and  
vertical lines onto unstable arcs.
The \emph{corners} and the \emph{sides} of the product box 
are the images of the corners and the sides of the square $[0,1]\times[0,1]$.
\end{df}

\begin{teo}
\label{boxprod}
If $f$ is a regular cw-expansive surface homeomorphism without arbitrarily small doubly asymptotic sectors then for all $x\in S$ 
a neighborhood of $x$ is covered by at least 4 product boxes with corner at $x$ and intersecting in the sides.
\end{teo}

\begin{proof}
Given $x\in S$ consider an open disc $D$ around $x$. 
Let $\alpha$ be a stable arc from $x$ to $\partial D$. 
Consider an unstable arc $\beta$ from $x$ to $\partial D$ such that in some component $U$ 
of $D\setminus (\alpha\cup \beta)$ there is 
neither a stable nor unstable arc from $x$ to $\partial D$ contained in $U$. 
By Proposition \ref{inestsepara} for each $y\in \alpha$ there is an unstable arc $u(y)\subset U$ 
from $y$ to $\partial U$. 
Also, for all $z\in\beta$ there is a stable arc $s(z)\subset U$ 
from $z$ to $\partial U$. 
Let us prove that there are $y_0\in\alpha$ and $z_0\in\beta$ such that 
$u(y_0)$ cuts $s(z_0)$. See Figure \ref{figProdLoc}.
\begin{figure}[ht]
 \begin{center}
  \includegraphics{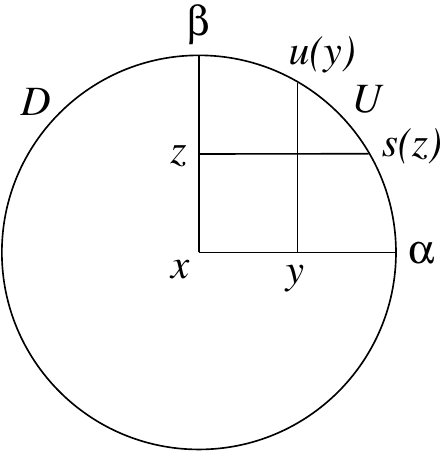}
  \caption{Local product structure in a sector of $x$.}
  \label{figProdLoc}
 \end{center}
 \end{figure}
If this is not the case we can take $y_n\in \alpha$ and $z_n\in\beta$ both converging to $x$ 
and such that $u(y_n)\cap s(z_n)=\emptyset$. 
Then we have a contradiction by taking limit in the Hausdorff metric of the continua 
$u(y_n)$ (or $S(z_n)$), because we obtain an unstable arc contained in $U$ from $x$ to $\partial D$. 
Let $V$ be the rectangle limited by the four 
curves $\alpha,\beta,u_{y_0}$ and $s_{z_0}$. 
Applying again Proposition \ref{inestsepara} in $V$ we have that the closure of 
$V$ is a product box. 
Now applying this construction a finite number of times we cover a neighborhood of $x$.
There will be at least 4 product boxes because there are no doubly asymptotic sectors. 
\end{proof}

\subsection{Expansive homeomorphisms}

\begin{teo}
\label{TeoEquivExp}
If $f\colon S\to S$ is a homeomorphism of a compact surface $S$ then
the following statements are equivalent:
 \begin{enumerate}
  \item $f$ is expansive
  \item $f$ is cw-expansive without topological doubly asymptotic sectors.
  \item $f$ is a regular cw-expansive homeomorphism without doubly asymptotic sectors,
  \end{enumerate}
\end{teo}

\begin{proof}
($1\to 2$).
Expansive homeomorphisms of compact metric spaces are always cw-expansive. 
Also, expansive homeomorphisms of compact surfaces are known to be conjugate to pseudo-Anosov diffeomorphisms. 
Then, it is easy to see that there are no topological doubly asymptotic sectors.

($2\to 3$). By Proposition \ref{cwNotbimpreg} we have that $f$ is a regular cw-expansive homeomorphism. 
It only rest to note that doubly asymptotic sectors are topological doubly asymptotic sectors.

($3\to 1$). 
By Theorem \ref{boxprod} we can see that there is $\alpha>0$ such that $W^s_\alpha(x)\cap W^u_\alpha(x)=\{x\}$ for all $x\in S$. 
Therefore, $\alpha$ is an expansive constant.
\end{proof}

\section{Anomalous expansive homeomorphisms}
\label{secAno}
In this section we will show that there are cw-expansive surface homeomorphisms and expansive homeomorphisms of three manifolds with a point 
whose local stable set is not locally connected.
Both examples are based on the construction of an irregular saddle in the plane.
\subsection{Irregular saddle}
\label{secIrregular}
Now we will construct a plane homeomorphism with a fixed point whose local stable set is not locally connected.

\begin{figure}[ht]
\center{\includegraphics{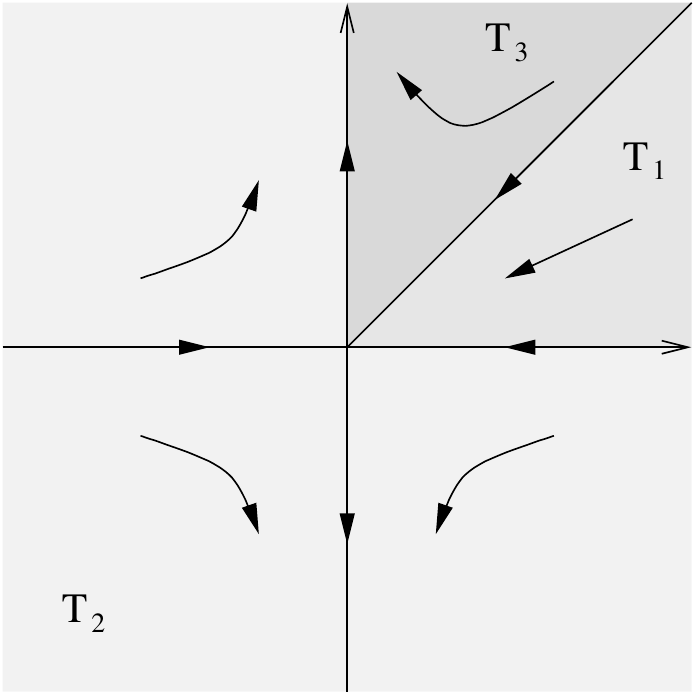}}
\caption{The action of the piece-wise linear transformation $T$.}
\label{sillaLoca}
\end{figure}

\begin{figure}[ht]
\center{\includegraphics{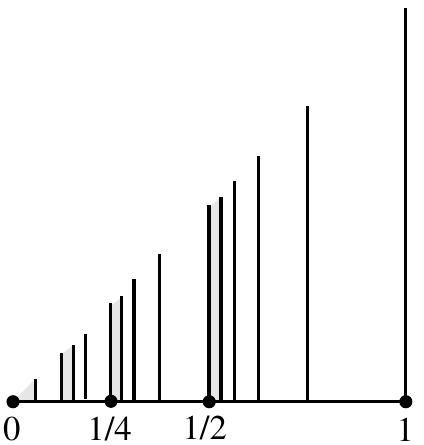}}
\caption{The set $E$.}
\label{conjuntoE}
\end{figure}

Let $T_i\colon \R^2\to\R^2$, for $i=1,2,3$, be the linear transformations defined by 
 $T_1(x,y)=(\frac x2,\frac y2)$, $T_2(x,y)=(\frac x2,2y)$, $T_3(1,1)=(\frac12,\frac12)$, $T_3(0,1)=(0,2)$.
Define the piece-wise linear transformation $T\colon \R^2\to\R^2$ as 
\[
 T(x,y)=\left\{ 
 \begin{array}{l}
    T_1(x,y)\hbox{ if } x\geq y\geq 0,\\
    T_2(x,y)\hbox{ if } x\leq 0 \hbox{ or } y\leq 0,\\
    T_3(x,y)\hbox{ if } y\geq x\geq 0.
 \end{array}
 \right.
\]
In Figure \ref{sillaLoca} we illustrate the definition of $T$.
Define the sets 
$C(a)=\{(a,y)\in\R^2:0\leq y\leq a\}$ for $a>0$, 
$D_1=\cup_{i=1}^\infty C(\frac12+\frac1{2^i})$, 
$D_{n+1}=T_1(D_n)$ for all $n\geq 1$. 
Finally consider the non-locally connected continuum $E=D\cup([0,1]\times\{0\})$ shown in Figure \ref{conjuntoE}. 

Consider the continuous function $\rho\colon \R^2\to\R$ defined by 
$$\rho(p)=\dist(p,E)=\min\{\dist(p,q):q\in E\}$$
and the vector field $X\colon\R^2\to\R^2$ defined as 
$$X(p)=(0,\rho(p)).$$ It is easy to see that $X$ has unique solutions, therefore, we can consider the flow $\phi\colon\R\times\R^2\to\R^2$ induced by $X$. 

Let $f\colon\R^2\to\R^2$ be the homeomorphism 
$$f=\phi_1\circ T,$$ 
where $\phi_1$ is the time-one homeomorphism associated to the vector field $X$.

\begin{prop}
 The homeomorphism $f$ preserves the vertical foliation on $\R^2$.
\end{prop}

\begin{proof}
 It follows because $\phi_t$ and $T$ preserves the vertical foliation.
\end{proof}

Consider the region $$R_1=\{(x,y)\in[0,1]\times[0,1]:x\geq y\}.$$ 

\begin{lem}
\label{lemaR1}
 For all $p\in R_1$ it holds that $\rho(T(p))=\frac 12 \rho(p)$ and $\phi_t(T(p))=T(\phi_t(p))$ if $\phi_t(p)\in R_1$ and $t\geq 0$.
\end{lem}

\begin{proof}
 By the definition of $T$ we have that $T(p)=T_1(p)=\frac12p$ for all $p\in R_1$. 
 Given $p\in R_1$ consider $q\in E$ such that $\rho(p)=\dist(p,q)$. 
 Then $\rho(T(p))=\dist(T(p),T(q))$ and $\rho(T(p))=\frac 12 \rho(p)$.
  
 Consider $t\geq 0$ such that $\phi_t(p)\in R_1$. Since $X$ is a vertical vector field we have that $\phi_{[0,t]}(p)\subset R_1$. 
 For $s\in(0,t)$, if $q=\phi_s(p)$ then 
 \[
  X(T(q))=(0,\rho(T(q)))=\left(0,\frac12\rho(q)\right)=d_qT(X(q)).
 \]
 Therefore, $\phi_s(T(p))=T(\phi_s(p))$ for $s\in(0,t)$ and consequently for $s=t$.
\end{proof}

Define the stable set of the origin as usual by $$W^s_f(0)=\{p\in\R^2:\lim_{n\to+\infty}\|f^n(p)\|=0\}.$$

\begin{prop}
 For the homeomorphism $f\colon\R^2\to\R^2$ defined above it holds that $$W^s_f(0)\cap[0,1]\times[0,1]=E.$$
\end{prop}

\begin{proof} 
First notice that $E\subset W^s_f(0)$ because for all $p\in E$ and $t\in\R$ we have that 
$\phi_t(p)=p$, and $T(p)=\frac12p$. Then $f(p)=\frac12p$ for all $p\in E$. 

Now take a point $p\in [0,1]\times[0,1]$. 
For $p\notin R_1$ it is easy to see that $f^n(p)\to\infty$ as $n\to+\infty$.
Assume that $p\in R_1\setminus E$. 
We will show that $p\notin W^s_f(0)$. 
It is sufficient to show that for some $n>0$ the point $f^n(p)$ is not in $R_1$.
By contradiction, assume that $f^n(p)\in R_1$ for all $n\geq 0$. 
Then, by Lemma \ref{lemaR1} we know that 
$$f^n(p)=(\phi_1\circ T)^n(p)=T^n(\phi_n(p)).$$
Notice that $\phi_1^n=\phi_n$. 
Then, it only rests to prove that $\phi_n(p)\notin R_1$ for some $n>0$. 
But this is easy because the velocity of $\phi_t(p)$ is $\rho(\phi_t(p))$ and 
this velocity increases with $t$. 
\end{proof}

\begin{rmk}[Three-dimensional irregular saddle]
\label{3dIrrSaddle}
 A three dimensional version of this example can be given as follows. 
 Consider the homeomorphism of $\R^3$ defined as $(x,y,z)\mapsto(f(x,y),2z)$. 
 It is easy to see that the local stable set (for the new map) of the origin is $E\times \{0\}$.
\end{rmk}

\subsection{Cw-expansive surface homeomorphisms}
\label{AnomCwexp}
The purpose of this section is to show that there are cw-expansive surface 
homeomorphisms with points whose local stable set is not locally connected.
The construction of this example is based on the construction  
of a quasi-Anosov diffeomorphism given in \cite{FR} using our irregular saddle.
It starts with a well know diffeomorphism that we will separate for future reference.

\subsubsection{Derived from Anosov}
\label{secDAnormal}
\begin{figure}[ht]
\begin{center}  
  \includegraphics{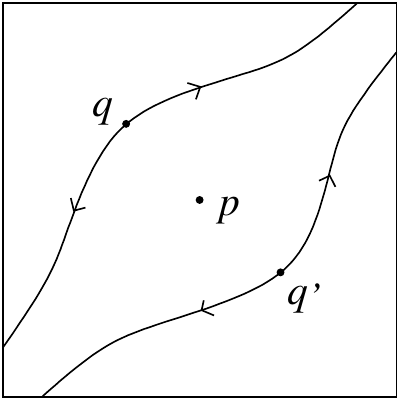}
  \caption{The derived from Anosov diffeomorphism on the two dimensional torus.}
  \label{figDA}
\end{center}
\end{figure}

In this section we will recall some properties of what is known as 
a derived from Anosov diffeomorphisms. 
The interested reader should consult \cite{Robinson} Section 8.8 for a construction of such a map 
and detailed proofs of its properties. 
A derived from Anosov is a $C^\infty$ diffeomorphism $f\colon T^2\to T^2$ 
of the torus such that: it satisfies Smale's axiom A and its non-wandering set consists 
of an expanding attractor and a repelling fixed point $p\in T^2$ as in Figure \ref{figDA}.

We will also assume that there is a local chart $\varphi\colon D\to T^2$, 
defined on the disc $D=\{x\in\R^2:\|x\|\leq 2\}$,
such that
\begin{enumerate}
\item $\varphi(0)=p$,
\item the pull-back of the stable foliation by $\varphi$
is the vertical foliation on $D$ and
\item $\varphi^{-1}\circ f\circ \varphi(x)=4x$
for all $x\in D$.
\end{enumerate}
More details of the construction can be found in \cite{FR}. 

\subsubsection{Irregular derived from Anosov}
\label{secIDA}
\begin{figure}[ht]
\begin{center}  
  \includegraphics{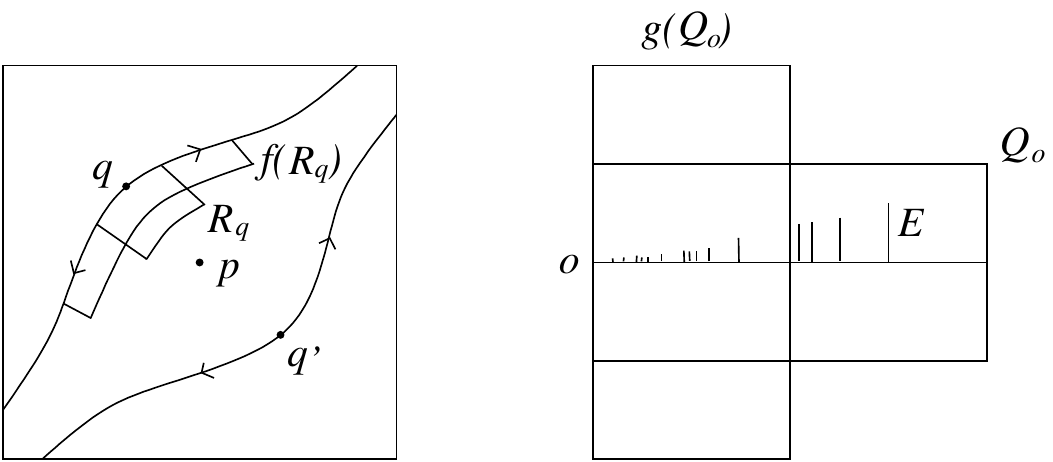}
  \caption{Topological rectangles on the derived from Anosov (left) and on the irregular saddle (right).}
  \label{figDARect}
\end{center}
\end{figure}

\begin{figure}[ht]
\begin{center}  
  \includegraphics{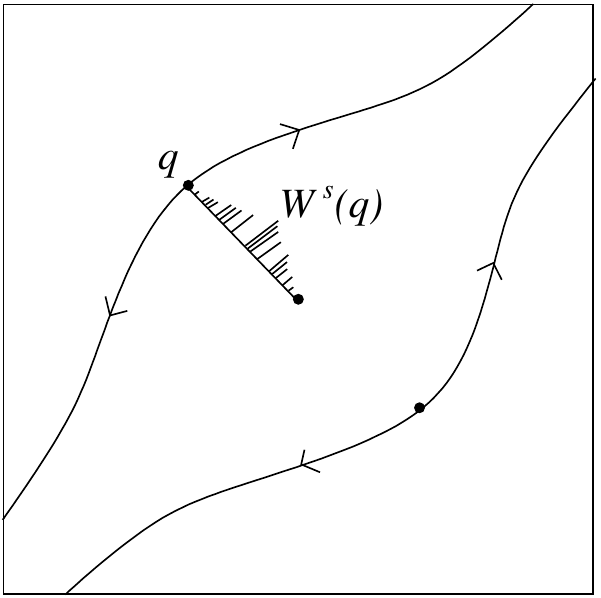}
  \caption{Derived from Anosov with an irregular saddle at the fixed point $q$.}
  \label{figAnoDA}
\end{center}
\end{figure}

Let $q$ and $q'$ be the two hyperbolic fixed points of the derived from Anosov diffeomorphism shown in Figure \ref{figDA}. 
Consider a topological rectangle $R_q$ covering a half-neighborhood of $q$ as in Figure \ref{figDARect}.
Consider the homeomorphism with an irregular saddle fixed point defined in Section \ref{secIrregular}. 
Call this homeomorphism $g$ (to avoid confusion with the derived from Anosov $f$). 
Denote by $o$ its fixed point (the origin of $\R^2$) and take a rectangle $Q_o$, similar to $R_p$, as in Figure \ref{figDARect}.
Now we can \emph{replace} $R_q$ with $Q_o$ and define what we call a 
\emph{derived from Anosov with an irregular saddle} as in Figure \ref{figAnoDA}. 

\subsubsection{Anomalous cw-expansive surface homeomorphism}
In this section we will construct a cw-expansive surface homeomorphism with a point whose local stable set 
is connected but it is not locally connected.
Consider $S_1$ and $S_2$ two disjoint copies of the torus $\R^2/\Z^2$.
Let $f_i\colon S_i\to S_i$, $i=1,2$, be two homeomorphisms such that:
\begin{itemize}
\item $f_1$ is the derived from Anosov with an irregular saddle from Section \ref{secIDA}, denote by  
$p_1\in S_1$ the source fixed point of $f_1$,
\item $f_2$ is the inverse of the derived from Anosov (the usual one as in Section \ref{secDAnormal}) with a sink fixed point at $p_2\in S_2$.
\end{itemize}
Consider local charts $\varphi_i\colon D_2\to S_i$, $i=1,2$, where $D_2$ is the compact 
disk
$$D_2=\{x\in\R^2:\|x\|\leq 2\},$$
such that:
\begin{enumerate}
 \item $\varphi_i(0)=p_i$,
\item the pull-back of the unstable foliation by $\varphi_2$
is the vertical foliation on $D_2$ and
\item $\varphi_1^{-1}\circ f^{-1}_1\circ \varphi_1(x)=\varphi_2^{-1}\circ f_2\circ \varphi_2(x)=x/4$
for all $x\in D$.
\end{enumerate}

\begin{figure}[ht]
\begin{center}  
  \includegraphics{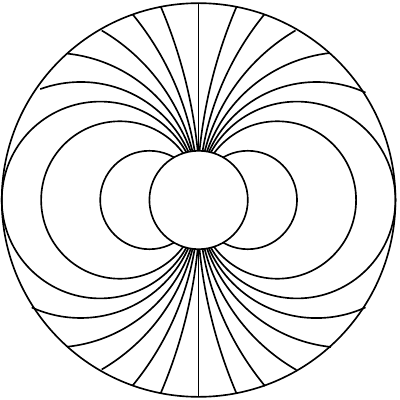}
  \caption{Unstable foliation of $f_2$ on the annulus, in the local chart $\varphi_2\circ\psi$.}
  \label{figFols}
\end{center}
\end{figure}

Consider the open disk $$D_{1/2}=\{x\in\R^2:\|x\|<1/2\}$$ 
and the compact annulus $$A=D_2\setminus D_{1/2}.$$
Define $\psi\colon A\to A$ as the inversion $\psi(x)=x/\|x\|^2$.
The pull-back of the unstable foliation on $S_2$ by $\varphi_2\circ\psi$ on the annulus $A$ 
is shown in Figure \ref{figFols}.
On the disjoint union $S_3=[S_1 \setminus \varphi_1(D_{1/2})]\cup [S_1\setminus \varphi_2(D_{1/2})]$
consider the equivalence relation generated by
\[
 \varphi_1(x)\sim \varphi_2\circ\psi (x)
\]
for all $x\in A$. Denote by $[x]$ the equivalence class of $x$.
The surface $S=S_3/\sim$ is the genus two surface if equipped with the quotient topology.
Consider the homeomorphism $f\colon S\to S$ defined by
\[
 f([x])=\left\{
\begin{array}{ll}
\left[f_1(x)\right] &\hbox{ if } x\in S_1 \setminus \varphi_1(D_{1/2})\\
\left[f_2(x)\right] &\hbox{ if } x\in S_2 \setminus \varphi_2(D_2)\\
\end{array}
\right.
\]
We will show that there are $C^0$ perturbations of $f$ that are cw-expansive 
and has a stable continuum that is not locally connected.

\begin{teo}
There are cw-expansive homeomorphisms of the genus two surface having stable continua that are not 
locally connected.
\end{teo}

\begin{proof}
Define $A_S=[\phi_1(A)]$ the annulus on $S$ corresponding to $A$. 
We will perturb the homeomorphism $f$ defined above on the annulus $A_S$. 
First note that the non-wandering set of $f$ is expansive and isolated. 
Also note that for every wandering point $x\in S$ there is $n\in\Z$ such that 
$f^n(x)\in A_S$. 
Therefore, it is sufficient to prove that there is a homeomorphism $g\colon S\to S$ such that 
$f|A_S=g|A_S$ and 
there is $\delta>0$ such that for each $x\in A_S$ the intersection 
$W^s_\delta(x)\cap W^u_\delta(x)$ is zero-dimensional, that is, it contains no non-trivial continua. 
In Figure \ref{figFols} we have the picture of the unstable foliation on $A_S$ (or in local charts). 
The problem is that the stable sets do not make a foliation, this is because there is an irregular saddle. 
Therefore, it is convenient to consider the stable partition, i.e., the partition defined by 
the equivalence relation of being positively asymptotic. 
This partition is illustrated in Figure \ref{figStPart}.
\begin{figure}[ht]
\begin{center}  
  \includegraphics{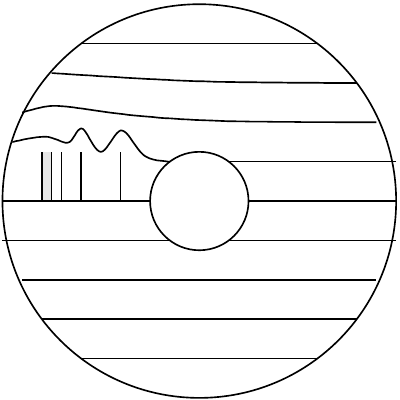}
  \caption{Stable partition on the annulus $A_S$.}
  \label{figStPart}
\end{center}
\end{figure}
We know that the unstable leaves are circle arcs, as in Figure \ref{figFols}. 
Therefore it is sufficient to consider a $C^0$ perturbation $g$ of $f$ supported on $A_S$, 
such that the stable partition of $g$ in the annulus contains no circle arc, in local charts. 
By the previous comments it implies that $g$ is cw-expansive because 
the intersection of local stable sets and local unstable sets are at most countable. 
Since $g$ coincides with $f$ outside $A_S$, we have that $g$ has an irregular saddle with 
non-locally connected stable set. This finishes the proof.
\end{proof}

\subsection{Expansive homeomorphisms on three-manifolds}

For expansive surface homeomorphisms it is known that every local stable set is locally connected, see \cites{Hi,L}.
The same is true for expansive diffeomorphisms of three-manifolds without wandering points as 
proved in \cite{Vi2002}.
In this section we will construct an expansive homeomorphism 
on a three-dimensional manifold with a 
point whose local stable set is not locally connected. 
The non-wandering set consists on an attractor and a repeller.

\begin{teo}
There is an expansive homeomorphism of a three-dimensional manifold with a point whose local stable set is connected but not locally connected.
\end{teo}

\begin{proof}
Let $f\colon M\to M$ be the quasi-Anosov diffeomorphism (that is not Anosov) constructed in \cite{FR}. 
Its non-wandering set contains an expanding attractor $\Omega_a$ and a contracting repeller $\Omega_r$. 
It is constructed by surgery on the wandering sets of two copies of three-dimensional derived from Anosov diffeomorphisms, 
\begin{figure}[ht]
\begin{center}  
  \includegraphics{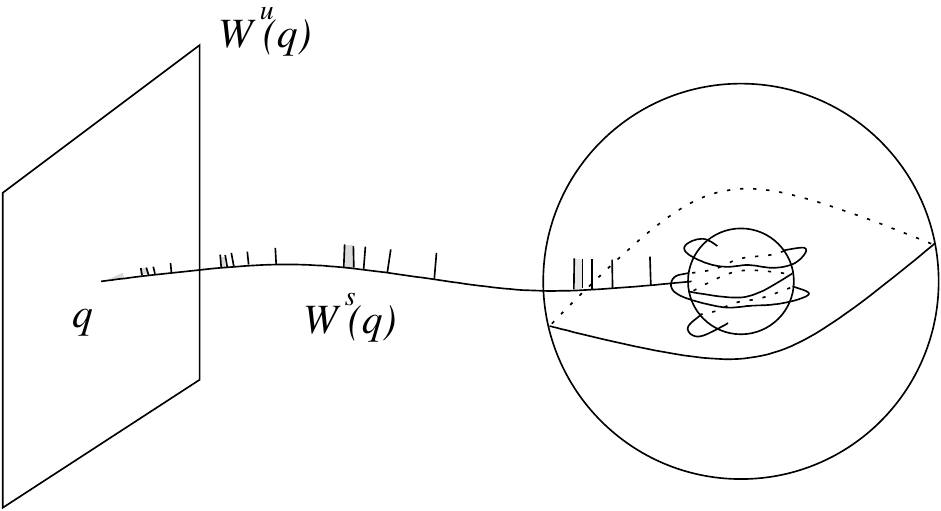}
  \caption{One the left there is a three-dimensional irregular saddle from Remark \ref{3dIrrSaddle}. 
  On the right there is a \emph{wormhole} with its \emph{twisted dipole foliation} from \cite{FR}. 
  The warm hole connects both derived from Anosov systems.}
  \label{figQAirr}
\end{center}
\end{figure}
therefore, there is a fixed point $q\in \Omega_a$ whose local unstable set separates a neighborhood $U$ of $q$ and 
one component of $U\setminus W_{loc}^u(q)$ is contained in the wandering set of $f$. 
On this component, \emph{replace} the diffeomorphism $f$ with 
a three-dimensional irregular saddle 
from Remark \ref{3dIrrSaddle}. 
It is easy to see that, eventually perturbing the obtained homeomorphism, we have that stable and unstable sets meets locally in at most one point (the hard part is done in \cite{FR}). 
In Figure \ref{figQAirr} a sketch is given. By definition, this implies that the homeomorphism is expansive. 
The point $q$ has a local stable set that is connected but not locally connected.
\end{proof}

\section{Stably expansive diffeomorphisms}
\label{secSta}
In this section we extend the definition of quasi-Anosov diffeomorphism 
by requiring some control in the order of a tangency of stable and unstable manifolds. 
We relate this concept with robust $N$-expansiveness.

\subsection{Omega-expansivity}

Let $M$ be a smooth compact surface without boundary. 
Given a diffeomorphisms $f\colon M\to M$ define 
$\per(f)$ as the set of periodic points of $f$ and the non-wandering set $\Omega(f)$ as the set of those $x\in S$ satisfying: 
for all $\epsilon>0$ there is $n\geq 1$ such that $B_\epsilon(x)\cap f^n(B_\epsilon(x))\neq\emptyset$.
Recall that $f$ satisfies Smale's \emph{Axiom A} 
if $\clos(\per(f))=\Omega(f)$ and $\Omega(f)$ is hyperbolic. 
Recall that $\Lambda\subset M$ is \emph{hyperbolic} 
if it is compact, invariant and
the tangent bundle over $\Lambda$ 
splits as $T_{\Lambda}M=E^s\oplus E^u$ 
the sum of two sub-bundles invariant by $df$
and there are $c>0$ and $\lambda\in (0,1)$ such that:
\begin{enumerate}
 \item if $v\in E^s$ then $\|df^n(v)\|\leq c\lambda^n \|v\|$ for all $n\geq 0$ and 
 \item if $v\in E^u$ then $\|df^n(v)\|\leq c\lambda^n \|v\|$ for all $n\leq 0$.
\end{enumerate}

\begin{df}
 We say that $f$ is $\Omega$-\emph{expansive} if $f\colon \Omega(f)\to\Omega(f)$ is expansive. 
 We say that $f$ is \emph{robustly} $\Omega$-\emph{expansive} if there is an 
 open $C^1$-neighborhood $U$ of $f$ such that every $g\in U$ is $\Omega$-expansive.
\end{df}

\begin{df}
 A $C^1$ diffeomorphism $f\colon M\to M$ is $\Omega$-stable if there is 
 a $C^1$ neighborhood $U$ of $f$ such that 
 for all $g\in U$ there is a homeomorphism $h\colon \Omega(f)\to\Omega(g)$ 
 such that $h\circ f=g\circ h$. 
\end{df}

Recall that $f$ is a \emph{star diffeomorphism} if there is a $C^1$ neighborhood $U$ of $f$ such that 
every periodic point of every $g\in U$ is hyperbolic.
If $f$ satisfies the axiom A then $\Omega(f)$ decomposes in a finite 
disjoint union basic sets $\Omega(f)=\Lambda_1\cup\dots\cup\Lambda_l$. 
A collection $\Lambda_{i_1},\dots,\Lambda_{i_k}$ is called a \emph{cycle} 
if there exist points $a_j\notin\Omega(f)$, for $j=1,\dots,k$, 
such that $\alpha(a_j)\subset \Lambda_{i_j}$ 
and $\omega(a_j)\subset \Lambda_{i_{j+1}}$ (with $k+1\equiv 1$). 
We say that $f$ \emph{has not cycles} 
if there are not cycles among the basic sets of $\Omega(f)$. 
See for example \cite{Robinson} for the definition of basic set and more on this subject.

\begin{teo}
 The following statements are equivalent:
 \begin{enumerate}
  \item $f$ satisfies axiom A and has not cycles,
  \item $f$ is $\Omega$-stable,
  \item $f$ is a star diffeomorphism.
 \end{enumerate}
\end{teo}

\begin{proof}
 ($1\Rightarrow 2$). It was proved by Smale in 1970 \cite{Sm70}. 
 
 ($2\Rightarrow 3$). It was proved by Franks in 1971 \cite{Fr71}. 

 ($3\Rightarrow 1$). It was proved in 1992 by Aoki \cite{Ao92} and Hayashi \cite{Ha}.
\end{proof}


\begin{teo}
 A $C^1$ diffeomorphism is robustly $\Omega$-expansive if and only if it is $\Omega$-stable.
\end{teo}

\begin{proof}
 If $f$ is $\Omega$-stable then $f$ satisfies Smale's axiom A. 
 Therefore $\Omega(f)$ is hyperbolic and consequently $f\colon\Omega(f)\to\Omega(f)$ is expansive. 
 Since $f$ is $\Omega$-stable we have that $f$ is robustly $\Omega$-expansive. 
 
 In order to prove the converse, suppose that $f$ is robustly $\Omega$-expansive. 
 Using Franks' Lemma it is easy to see that $f$ is a star diffeomorphism.
\end{proof}

\subsection{Stable N-expansive surface diffeomorphisms}
In this section we will consider a diffeomorphism $f$ of a $C^\infty$ compact surface $S$.
The \emph{stable set} of $x\in S$ is
\[
 W^s_f(x)=\{y\in S:\lim_{n\to+\infty}\dist(f^n(x),f^n(y))= 0\}.
\]
The \emph{unstable set} is defined by $W^u_f(x)=W^s_{f^{-1}}(x)$.
Assume that $f$ is $\Omega$-stable, $E^s, E^u$ are one-dimensional
and define $I=[-1,1]$.
Let us recall the following fundamental result for future reference.

\begin{teo}[Stable manifold theorem]
\label{teoSMT}
 Let $\Lambda\subset S$ be a hyperbolic set of a 
 $C^r$ diffeomorphisms $f$ of a compact surface $S$. 
 Then, for all $x\in\Lambda$, $W^s_f(x)$ is an injectively immersed $C^r$ submanifold. 
 Also the map $x\mapsto W^s_f(x)$ is continuous: 
 there is a continuous function $\Phi\colon\Lambda\to\Emb^r(I,S)$ 
 such that for each
 $x\in\Lambda$ it holds that the image of $\Phi(x)$ is a neighborhood of $x$ in $W^s_f(x)$. 
 Finally, these stable manifolds also depend continuously on the diffeomorphisms $f$, 
 in the sense that nearby diffeomorphisms yield nearby mappings $\Phi$ as above.
 \end{teo}
\begin{proof}
 See  \cite{PaTa} Appendix 1.
\end{proof}

\begin{df}
A $C^r$, $r\geq 1$, diffeomorphisms $f\colon S\to S$ is \emph{$Q^r$-Anosov} 
if it is $\Omega$-stable in the $C^r$ topology and for all $x\in S$ 
there are $\delta_1,\delta_2>0$, a $C^r$ coordinate chart 
$\varphi\colon U\subset S\to [-\delta_1,\delta_1]\times[\delta_2,\delta_2]$ 
such that $\varphi(x)=(0,0)$ and 
two $C^r$ functions $g^s,g^u\colon [-\delta_1,\delta_1]\to[\delta_2,\delta_2]$ such that 
the graph of $g^s$ and $g^u$ are the local 
expressions of the local stable and the local unstable manifold of $x$, respectively, 
and the degree $r$ Taylor polynomials of $g^s$ and $g^u$ at $0$ are different. 
If the polynomials coincide we say that there is an $r$-\emph{tangency} at the intersection point.
\end{df}

\begin{rmk}
 For $r=1$ we have that $Q^1$-Anosov is quasi-Anosov, and in fact, given that $S$ is two-dimensional, it is Anosov. 
 For $r=2$ we are requiring that if there is a tangency of a stable and an unstable manifold it is a quadratic one.
\end{rmk}

\begin{lem}
\label{lemPertOmega}
 If $f$ is $\Omega$-stable then for all $\epsilon>0$ there are $m\geq 0$  
 and a $C^1$ neighborhood ${\cal U}$ of $f$ such that 
 if $|n|\geq m$ then
 $g^n(x)\in B_\epsilon(\Omega(g))$
 for all $x\in M$ and $g\in {\cal U}$.
\end{lem}

\begin{proof}
 By contradiction, take $\epsilon>0$, $g_k\to f$, $x_k\in M$
 and $n_k\to\infty$ such that 
 for all $k\in\N$, $\dist(g_k^i(x_k,\Omega(g_k))\geq \epsilon$ if $|i|\leq n_k$. 
 By Theorem 8.3 in \cite{Sh87} and the $\Omega$-stability of $f$ we know that $\Omega(g_k)\to \Omega(f)$ in the Hausdorff metric.
 Therefore, if $x_k\to x$ then $\dist(f^i(x),\Omega(f))\geq\epsilon$ for all $i\in\Z$. 
 But this is a contradiction because $\omega_f(x)\subset \Omega(f)$.
\end{proof}

\begin{teo}
 In the $C^r$ topology 
 the set of $Q^r$-Anosov diffeomorphisms is an open set. 
\end{teo}

\begin{proof}
We know that the set of $\Omega$-stable diffeomorphisms is an open set. 
Let $g_k$ be a sequence of $\Omega$-stable $C^r$-diffeomorphisms 
converging to the $C^r$ $\Omega$-stable diffeomorphism $f$. 
Assume that $g_k$ is not $Q^r$-Anosov for all $k\geq 0$. 
In order to finish the proof is it sufficient to show that $f$ is not $Q^r$-Anosov. 
Since $g_k$ is $\Omega$-stable but it is not $Q^r$-Anosov, there is $x_k\in S$ with an $r$-tangency. 
Assume that $x_k\to x$. By Theorem \ref{teoSMT} and Lemma \ref{lemPertOmega} we have that there is an $r$-tangency at 
$x$. Therefore $f$ is not $Q^r$-Anosov. 
Consequently, the set of $Q^r$-Anosov $C^r$-diffeomorphisms is an open set in the $C^r$ topology.
\end{proof}

Let us recall a definition from \cite{Mo12}.

\begin{df}
 Given $N\in\Z^+$ we say that a homeomorphism $f\colon S\to S$ is $N$-\emph{expansive} if 
 there is $\delta>0$ such that if $\diam(f^i(A))\leq\delta$ for all $i\in\Z$ then 
 $A$ has at most $N$ points. 
\end{df}

\begin{rmk}
 A homeomorphism is 1-expansive if and only if it is expansive. 
 Also, every $N$-expansive homeomorphism is cw-expansive.
\end{rmk}

\begin{df}
  We say that a $C^r$ diffeomorphism $f$ is $C^r$-\emph{robustly} $N$-expansive if 
 there is a $C^r$ neighborhood of $f$ such that every diffeomorphism in this neighborhood is $N$-expansive.
\end{df}

\begin{rmk}
For $r=N=1$ Mañé \cite{Ma75} proved that a diffeomorphism is robustly expansive if and only if it is quasi-Anosov. 
\end{rmk}

\begin{lem}
 If $g\colon \R\to\R$ is a $C^r$ functions with $r+1$ roots in the interval 
 $[a,b]\subset\R$ then $g^{(n)}$ has $r+1-n$ roots in $[a,b]$ for all $n=1,2,\dots,r$ where $g^{(n)}$ stands for the $n^{th}$ derivative of $g$.
\end{lem}

\begin{proof}
 It follows by induction in $n$ using 
 the Rolle's theorem.
\end{proof}

\begin{teo}
\label{teoRobNexp}
Every $Q^r$-Anosov $C^r$ diffeomorphism of a compact surface is 
 $C^r$-robustly $r$-expansive.
\end{teo}

\begin{proof}
It follows by definitions and the previous Lemma.
\end{proof}

\begin{rmk} Theorem 
\ref{teoRobNexp} can be applied 
in the example constructed in \cite{APV}. 
This is because in local charts the tangencies in the wandering set of this example 
are made by straight lines and circle arcs. 
Consequently, they are quadratic tangencies.
This example is robustly 2-expansive (and cw-expansive) in the $C^2$ topology.
\end{rmk}
\subsection{Robust sensitivity}

Let $f\colon M\to M$ be a $C^1$ diffeomorphism of an $n$-dimensional compact manifold $M$.

\begin{teo}
 Every $\Omega$-stable diffeomorphism without periodic sinks or 
 without periodic sources is $C^1$-robustly sensitive to initial conditions.
\end{teo}


\begin{proof}
Let $\delta>0$ be so that the diameter of the stable set and the unstable set of every periodic point 
is greater than $2\delta$. 
Let $U\subset M$ be an open set. 
Since $f$ is Axiom A we have that the union of stable and unstable sets of periodic points is dense in $M$. 
By hypothesis, there is $x\in U$ such that $x\in W^\sigma(p)$ with 
$p$ periodic, $\sigma=s$ or $u$ and $\dim(W^\sigma(p))<n$. 
Then, there is a submanifold $N\subset U$ transverse to $W^\sigma(p)$ at $x$ such that $\dim(N)\geq 1$.
Then, applying the $\lambda$-Lemma we have that there is $n\in\Z$ such that $\diam(f^n(N))>\delta$. 
Consequently, $\diam(f^n(U))>\delta$ and $f$ is sensitive to initial conditions.
The robustness is a consequence of the fact that the set of diffeomorphisms satisfying the hypothesis is open.
\end{proof}

\begin{rmk}
The derived from Anosov diffeomorphisms satisfies the hypothesis of the theorem. 
They are not expansive but has robust sensitivity to initial conditions.
\end{rmk}

\noindent Departamento de Matemática y Estadística del Litoral, Salto-Uruguay\\
Universidad de la República\\
E-mail: artigue@unorte.edu.uy
\end{document}